\documentclass[letterpaper, 10 pt, conference]{ieeeconf}  

\IEEEoverridecommandlockouts                              

\title{\LARGE \bf
Stochastic Error Bounds in Nonlinear Model Predictive Control with Gaussian Processes via Parameter-Varying Embeddings  
}
 
\author{Dimitrios S. Karachalios$^{1}$ and Hossam S. Abbas$^{1}$
\thanks{This work has received funding from the Deutsche Forschungsgemeinschaft (DFG, German Research Foundation)-Project No. 419290163}
\thanks{$^{1}$The authors are with the Faculty of Electrical Engineering in Medicine, University of Luebeck, Germany. {\tt\small email:~dimitrios.karachalios@uni-luebeck.de}}%
}

\newtheorem{remark}{Remark}[section]

\def\IR{{\mathbb R}}

\def\IZ{{\mathbb Z}}

\newcommand{\cP}{ {\cal P} }

\def\IR{{\mathbb R}}

\usepackage{graphicx}
\usepackage{amsmath}
\usepackage{amsfonts}
\usepackage{amssymb}
\usepackage{svg}
\usepackage{times} 
\usepackage[capitalise]{cleveref}
\usepackage{algpseudocode}
\usepackage{algorithm}
\usepackage{caption}

\DeclareMathOperator{\sinc}{sinc}
\newtheorem{problem}{Problem}
\begin{document}
\maketitle
\thispagestyle{empty}
\pagestyle{empty}
\begin{abstract}
This study utilized the Gaussian Processes (GPs) regression framework to establish stochastic error bounds between the actual and predicted state evolution of nonlinear systems. These systems are embedded in the linear parameter-varying (LPV) formulation and controlled using model predictive control (MPC). Our main focus is quantifying the uncertainty of the LPVMPC framework's forward error resulting from scheduling signal estimation mismatch. We compared our stochastic approach with a recent deterministic approach and observed improvements in conservatism and robustness. To validate our analysis and method, we solved the regulator problem of an unbalanced disk.
\end{abstract}
\section{INTRODUCTION} 
To properly handle the inaccuracy from modeling mismatches between the predicted model and the real plant's response in model predictive control, research has been performed to create robust-tube control approaches for both linear and nonlinear systems in \cite{Mayne2011tubes,MAYNE2005219,Cannon2011,LANGSON2004125,Schildbach2013}. One of the many nonlinear control strategies that have drawn more attention is nonlinear model predictive control (NMPC), which converts the control task into an optimization problem considering input-state constraints. Since the MPC problem may be reformulated as an optimization problem of a quadratic manifold with an intrinsic unique optimal solution that can be efficiently found, embedding the nonlinear system in a linear parameter-varying (LPV) formulation is a great alternative to NMPC as the LPVMPC method solves the nonlinear control problem at time $k$ after only solving a single quadratic program (QP). Therefore, it can be handled in real-time, such as in autonomous driving tasks~\cite{Nezami2023}, using a variety of efficient algorithms.

One issue in using the LPV embedding to represent nonlinear systems is appropriately estimating the uncertainty quantification caused by the scheduling parameter $p$. The LPV operator realizes the nonlinear manifold via an adaptive tangential hyperplane, and the vector $p$ can absorb the nonlinear dependencies in an affine representation. The prediction of the scheduling parameter leads to an error between the actual and LPV model-predicted responses, which motivated this work.
%

Recently, there has been a lot of focus on the LPV models integrated within the MPC framework for handling the inherent uncertainty on $p$. One tactic for the LPV setting is to employ tube-based MPC, as proposed in~\cite{HanemaTubes}. This paper assumes the boundedness of future $p$ and proposes an anticipative tube MPC approach for LPV systems. Furthermore, a cascade MPC control architecture combining linear longitudinal and LPV lateral models was proposed in~\cite{HanemaTubes}, where anticipative tube MPC was proposed. This method inspired the description of an autonomous lane-keeping technique in~\cite{nezami2022robust}. Furthermore, a bilinear matrix inequality is used in~\cite{Abbas2016} to ensure the stability of the LPVMPC architecture. This approach, however, can be computationally costly and conservative. An offline technique was suggested for exploiting model uncertainty bounds as a function of the control inputs introduced in \cite{bujarbaruah2022robust}, where the MPC's tightening of constraints depends on decision variables or control inputs. 

Gaussian Processes (GPs) can be expensive because of the cubic complexity of the data set size. Nevertheless, substitutes, including the recently introduced Gaussian Processes Toolkit (GaPT) that does regression in real-time, have been showcased in \cite{Crocetti2023}. The GPs must be transformed into a linear state space form to scale linearly with data volume. Previous research dealing with disturbances \cite{Gruner2022}, careful control \cite{HeKaZe20}, and LPV \cite{Elkamel2022} has effectively used GPs. Furthermore, even while GPs employ a limited memory of trained data in a Bayesian conditional manner, the potential for active learning, when only hyper-parameter tuning is done, can further improve prediction in the sense of extrapolation. 

In \cref{sec:Prel}, we start with preliminaries and formal representation of the problem under consideration. In \cref{sec:MPCLPV-GPs}, we formulate explicitly the error dynamics that we want to predict together with the GPs regression framework. In \cref{sec:results}, we apply the derived new method to a classical control benchmark. Finally, in \cref{sec:conclusion}, we summarize our findings and provide the open challenges and future research directions.

\section{Preliminaries and Problem formulation}\label{sec:Prel}
\subsection{Definitions and assumptions}
We start with the discrete nonlinear dynamical system 
\begin{equation}\label{eq:sysnlc}
    \Sigma:x_{k+1}=f(x_k,u_k),
\end{equation}
of state dimension $n_{\mathrm{x}}$, and input dimension $n_u$. Considering the sampling time $t_s$, it holds $t_k=t_sk,~\forall k\in\IZ_+$, with $x_k=x(t_sk)$, and $x_0=x(0)$ the initial condition state vector. $f:\IR^{n_{\mathrm{x}}}\times\IR^{n_u}\rightarrow\IR^{n_{\mathrm{x}}}$ is a nonlinear operator. The nonlinear dynamical system in \eqref{eq:sysnlc} can be represented equivalently with a linear parameter-varying (LPV) formulation that will give rise to methods in an adaptive way.
\begin{equation}\label{sys:LPV}
    \Sigma:\left\{\begin{aligned}
        x_{k+1}&=A(p_k)x_k+Bu_k,\\
              p_k&=\rho(x_k),~x_0=x(0),\end{aligned}\right.
\end{equation}
An appropriate scheduling parameter vector $p$ with dimension $n_{p}$ should be introduced. The $p$ belongs in the scheduling set $\mathcal{P}:=\texttt{Co}\{p^{\nu_1},p^{\nu_2},\ldots,p^{\nu_{n_p}}\}$, where we denote with $p^{\nu_i}\in\IR^{n_p}$ the vertices of the set $\mathcal{P}$ based on given bounds of the scheduling $p$. Towards one further simplification, the remaining LPV system and through a pre-filter \cite{Apkarian1995} that will increase the state dimension can recast the scheduling dependence only to the linear matrix $A(p)$ allowing a non-dynamical-static matrix $B$. Consequently, the original nonlinear system \eqref{eq:sysnlc} can be embedded in the linear parameter varying (LPV) formulation  as in \eqref{sys:LPV}.

where the map $\rho:\IR^{n_{\mathrm{x}}}\rightarrow\IR^{n_{p}}$ is also given explicitly. In particular, the $\rho(\cdot)$ is a known nonlinear function of the state-$x$, which allows the embedding of \eqref{eq:sysnlc} in \eqref{sys:LPV}. Furthermore, the following parameterized matrix $A(p_k):\IR^{n_p}\rightarrow\IR^{n\times n}$ is also known and affine concerning the scheduling parameter $p$. In particular, the affine structure of the discrete operator $A(p_k)$ can be expressed as
\begin{equation}\label{eq:affine}
    A(p_k):=A_0+\sum_{l=1}^{n_p}p_k^{[l]}A_{l},
\end{equation}
where $p_k^{[l]}$ denotes the $l^{th}$-element of the vector $p$ and $A_l$ are constant matrices. Together with the input matrix $B\in\IR^{n_{\mathrm{x}}\times n_u}$, the discrete-time LPV system is well-defined. The following remark \ref{rem:stanassu} summarizes the appropriate general assumptions to solve the control task by considering LPV predictive models.
\begin{remark}[Standing assumptions]\label{rem:stanassu}
To reflect the generalization of the method, here are the minimal assumptions:
    \begin{itemize}
        \item Appropriately smoothness of the nonlinear operator $f$ has been assumed (i.e., higher-order differentiability and continuity).
        \item The scheduling parameter can be measured at each sampling time $k$ but remains unknown within the receding horizon of length N.
        \item The input matrix $B$ does not depend on the scheduling parameter. This can be relaxed easily through a pre-filter as in \cite{Apkarian1995}.
        \item No other disturbances or measurement noise has been assumed to affect the system $\Sigma$.
        \item The scheduling signal depends only on the state.
    \end{itemize}
\end{remark}
\subsection{Stabilization and model predictive control (MPC)} 
The dynamics around desired equilibrium operational points are unstable in many interesting control applications. Therefore, control strategies that stabilize the plant and drive the system to desired states that respect physical constraints are crucial. In the linear case, linear quadratic regulation (LQR) as a state feedback control strategy (i.e., $u=Kx$) efficiently provides stable closed-loop systems. In addition, LQR has been extended to handle LPV representations and provides robust state feedback that stabilizes the system for all possible parametrization of the scheduling $p\in\cP$ as in \cite{PANDEY20178624}. Thus, we propose splitting the control input into two parts; the first part will stabilize the nonlinear plant through a state-feedback auxiliary controller that can stabilize the LPV system over the whole range of $\mathcal{P}$ implemented in \cite{PANDEY20178624}, and the second part will handle the input-state constraints from an MPC control loop. In particular, the input design has the following structure:
\begin{equation}\label{eq:LQR}
    u_k=u_k^{\text{LQR}}+u_k^{\text{MPC}}=K x_k+u_k^{\text{MPC}}.
\end{equation}
By substituting \eqref{eq:LQR} to \eqref{sys:LPV}, we get
\begin{equation}
\footnotesize
    \begin{aligned}
        x_{k+1}&=A(p_k)x_k+Bu_k=A(p_k)x_k+B(K x_k+u_k^{\text{MPC}}),\\
         &=\underbrace{\left(A(p_k)+BK\right)}_{A_c(p_k)}x_k+Bu_k^{\text{MPC}}.\\
    \end{aligned}
\end{equation}
Thus, the closed-loop dynamics of applying the above auxiliary controller makes $A_c(p_k)$ stable and can be written as
\begin{equation}\label{sys:LPVLQR}
    \Sigma_{\text{c}}:\left\{\begin{aligned}
        x_{k+1}&=A_c(p_k)x_k+Bu_k,\\
              p_k&=\rho(x_k),
    \end{aligned}\right.
\end{equation}
with $A_c(p_k):=A(p_k)+BK$. We denote with subscript "c" the closed-loop operator, and the remaining controller $u_k ^{\text{MPC}}$ can be denoted again as $u_k$ without asserting any confusion.\\
\begin{problem}[Uncertainty quantification of the error]
   We want to establish stochastic bounds on the error between the true state response $x$ of the actual nonlinear system and the predicted state $\hat{x}$ within the prediction horizon of the MPC control strategy. The error is denoted as $e:=x-\hat{x}$. To address this, we will briefly explain the GP regression framework.
\end{problem} 
\subsection{Gaussian Processes (GPs)}
Gaussian Process regression is commonly used in machine learning due to its ability to measure uncertainty and incorporate prior knowledge about the data. It provides a non-parametric framework for model estimation, representing a probability distribution over functions \cite{RaWi05}. The following demonstration shows how GP regression can be used to learn $e=x-\hat{x}$. The continuous-time error $e(t)$ is modeled as a distribution over functions using a temporal, zero-mean prior $g_{\text{prior}}(t)\sim\mathcal{G}\mathcal{P}(0,k_{\text{prior}}(t,t'))$, where $k_{\text{prior}}(t,t')$ represents prior knowledge about the stochastic properties of the error. In this case, the squared exponential kernel is commonly chosen as the covariance function
\begin{equation}
    k_{\text{QP}}(t,t')=\sigma_{\text{QP}}^2 k_{\text{SE}}(t,t')k_{\text{P}}(t,t').
\end{equation}
Given the GP prior, a noisy measurement model for the error can be defined 
\begin{equation}
    \hat{e}_{j}:=g_{\text{prior}}(t_j)+\epsilon_j,
\end{equation}
where $\epsilon_j\sim \mathcal{N}(0,\Sigma^{\epsilon})$ represents i.i.d. measurement noise. Inferring the posterior GP distribution corresponds to solving the GP regression, which is done by conditioning on the measurements; see \cite{RaWi05} for details. We denote the learned posterior as $g_{\text{post}}(t)\sim \mathcal{G}\mathcal{P}(\bar{g}_{\text{post}},k_{\text{post}}(t,t'))$ and use it to predict the error as mentioned above. The GP's covariance functions usually have a small number of free hyper-parameters optimized beforehand by maximizing the marginal log-likelihood on a training set \cite{RaWi05}. 

After obtaining the hyper-parameters, we use Gaussian conditioning on the training set to ascertain the posterior GP distribution of the modeled function. The notation $g(t)\sim \mathcal{G}\mathcal{P}(\bar{g}(t),k^g_{\text{post}}(t,t'))$ represents the estimated posterior GP distribution. The posterior predictive distribution is then $p(g(t^*)|t,t*,y^g,\sigma^2)\sim \mathcal{N}(\bar{g}^*,\Sigma^{g^*})$, where the posterior mean $\bar{g}^*$ and covariance $\Sigma^{g^*}$ are provided by \cite{Klenske2016}. This posterior predictive distribution may estimate the function at a test input $t^*$ as
\begin{equation}
    \begin{aligned}
        &\bar{g}^*=k^T_{t^*t}K^{-1}y^g,\\
        &\Sigma^{g^{*}}=k_{t*t*}-k_{t^*t}K^{-1}k_{tt},
    \end{aligned}
    \end{equation}
where we use the following notation $k_{ab}=k^g(a,b)$ for simplicity.
In what comes next, we will use the GPs as in the following high-level information
\begin{equation*}
\footnotesize
    X_{t}\rightarrow\boxed{\textbf{GP}\text{(training)}}\rightarrow Y_{t}~~\&~X_{t^+}\rightarrow\boxed{\textbf{GP}\text{(trained)}}\rightarrow Y_{t^+}.
\end{equation*}
\section{Model predictive control with linear parameter-varying via Gaussian Processes (GPs)}\label{sec:MPCLPV-GPs}
\subsection{Error analysis within the prediction horizon}
The scope of the following analysis is to provide the error uncertainty quantification of the deviation between $\hat{x}_{i|k}$ (predicted state) and $x_{i|k}$ (true state) when a scheduling signal $\hat{p}_{i|k}$ has been used from a previous prediction within the prediction horizon. Once a scheduling parameter (signal) has been assumed, the solution of the QP offers the optimally designed input $u^\ast$ along with the predicted states $\hat{x}$. Implementing the input $u^\ast$ to the real system results in the true states $x$. These dynamical systems can be represented formally before and after the solution of the QP problem for $i=0,\ldots,N-1$ as  
    \begin{align}
    \footnotesize
\hat{x}_{i+1|k}&=A_c(\hat{p}_{i|k})\hat{x}_{i|k}+Bu_{i|k}^{\ast},~\text{(before QP solution)}\label{eq:LPVhat}\\
x_{i+1|k}&=A_c(p_{i|k})x_{i|k}+Bu_{i|k}^{\ast},~\text{(after QP solution)}\label{eq:LPVtrue}
    \end{align}
The error $\hat{e}_{i+1|k}\approx x_{i+1|k}-\hat{x}_{i+1|k}$ is predicted from the GPs and substituting $\hat{x}_{i+1|k}$ to \eqref{eq:LPVhat} remains
\begin{equation}\label{eq:stateImprov}
    x_{i+1|k}\approx\bar{x}_{i+1|k}=A_c(\hat{p}_{i|k})\hat{x}_{i|k}+Bu_{i|k}^{\ast}+\hat{e}_{i+1|k}.
\end{equation}
Learning the error in \eqref{eq:stateImprov} will make the estimated state closer to the true response, and subsequently, the inferred scheduling signal will be improved. In addition, as the error will be learned from a Bayesian framework as GPs, together with the likelihood estimation, higher statistical moments (i.e., variance) are offered. Consequently, the uncertainty quantification of the error can be determined. 

Quantifying the (parametric) uncertainty and prediction of the scheduling signal is at the heart of this study. The error source has to do with the fact that the scheduling estimation $\hat{p}_{i|k}$ is not necessarily the correct (true) one that the nonlinear plant will follow. 

In particular, by implementing the optimized control input $u_{i|k}^\ast$ to the real plant, the true response $x_{i|k}$ will deviate from the predicted one $\hat{x}_{i|k}$. Thus, we aim to provide bounds for this mismatch by introducing the error state at time $k$ and prediction $i$ (i.e., at $(i+k)^{\text{th}}$ real simulation time) between the true and predicted response that is defined as
\begin{equation}\label{eq:errorstate}
    e_{i|k}:=\underbrace{x_{i|k}}_{\text{true}}-\underbrace{\hat{x}_{i|k}}_{\text{predicted}},
\end{equation}
where $e_{0|k}=x_{0|k}-\hat{x}_{0|k}=x_0-x_0=0,~\forall k\in\IZ_+$. Once we get a solution $(\hat{x}_{i|k-1},u_{i|k-1}^\ast)$ at time $k-1$ from the optimization problem in \cref{eq:LPV_MPC}, we can implement the designed optimal input $u_{i|k-1}^\ast$ to the LPV that is equivalent with the original nonlinear and get the true responses. Therefore, we can compute explicitly the error $e_{i|k-1}$ that occurred at time step $k-1$ and within the prediction horizon for $i=1,\ldots,N$ that will consist of the training data for the GPs framework that will allow prediction of the error at the next time instant $k$ and before solving the MPC problem at time $k$.
\subsection{Error quantification with Gaussian Processes (GPs)}
At the heart of the proposed method lies the prediction of the forward error $e_{i+1|k}$ from past error sequences, i.e., $e_{i|k-1}$ with the GPs regression framework that subsequently will improve the estimation of the true response together with the uncertainty quantification. Such an effort can improve the robustness of the method against disturbances. We provide the following table consistent with the error \eqref{eq:errorstate}. The quadratic exponential kernel has been chosen for the GPs training phase in \cref{tab:errorTable} 
\vspace{5mm}
\begin{table}[h!]
\footnotesize
    \centering
    \begin{tabular}{c|ccc|c}
        $i\backslash k$ & $1$ & $2~\cdots$ & $k-1$ & $k$\\\hline\\[-2.5mm]
         $1$ & $e_{1|1}^{(n)}$ & $e_{1|2}^{(n)}~\cdots$  & $e_{1|k-1}^{(n)}$ & $\rightarrow\boxed{\textbf{GP}_1}\rightarrow e_{1|k}^{(n)}$ \\[2mm]
         $2$ & $e_{2|1}^{(n)}$ & $e_{2|2}^{(n)}~\cdots$ & $e_{1|k-1}^{(n)}$ & $\rightarrow\boxed{\textbf{GP}_2}\rightarrow e_{2|k}^{(n)}$\\[2mm]
         $3$ & $e_{3|1}^{(n)}$ & $e_{3|2}^{(n)}~\cdots$ & $e_{1|k-1}^{(n)}$ & $\rightarrow\boxed{\textbf{GP}_3}\rightarrow e_{3|k}^{(n)}$ \\[2mm]
         $\vdots$ & $\vdots$ & $\vdots~~~\ddots$ &  $\vdots$ & $\vdots$ \\[2mm]
         $N$ & $e_{N|1}^{(n)}$ & $e_{N|2}^{(n)}~\cdots$ & $e_{N|k-1}^{(n)}$ & $\rightarrow\boxed{\textbf{GP}_N}\rightarrow e_{N|k}^{(n)}$\\[2mm]\hline
         $N+1$ & $e_{N+1|1}^{(n)}$ & $e_{N+1|2}^{(n)}$ & $e_{N+1|k-1}^{(n)}$ & $\star$
    \end{tabular}
    \caption{At the time $k-1$, with prediction horizon length $N$ and for the error state $n=1,\ldots,n_{\textrm{x}}$, the prediction of the one step forward error $e_{i|k}^{(n)}$ with the GP framework is presented. The $\star$ indicates that at the time $k-1$, the measurement $e_{N+1|k}$ is not yet available, and neither prediction can be obtained.}
    \label{tab:errorTable}
\end{table}
where, the columns from $1$ to $k-1$ consist of the training data. In particular, the prediction mechanism that will be tackled with GPs for each error state is the following:
\begin{equation}\label{eq:GPlearnPredict}
    e_{i+1|k-1}^{(n)}\rightarrow\boxed{\begin{array}{cc}
         \textbf{GP}_{i}~\text{(training data)}  \\\hline
         \begin{array}{c|c}
                     X_i & Y_i\\\hline
             e_{i+1|1}^{(n)} & e_{i|2}^{(n)} \\
             e_{i+1|2}^{(n)} & e_{i|3}^{(n)} \\
             \vdots & \vdots \\
             e_{i+1|k-2}^{(n)} & e_{i|k-1}^{(n)} 
         \end{array}
    \end{array}}\rightarrow e_{i|k}^{(n)}~\text{(prediction)}
\end{equation}
To make the process in \eqref{eq:GPlearnPredict} handy, suppose we want to predict at time $k=2$ the forward error at $k=3$ and for the 1st point within the prediction horizon, i.e., $i=1$. Schematically, the prediction task is: $e_{2|2}\rightarrow\boxed{\textbf{GP}}\rightarrow e_{1|3}=?$. The GP training data consists of the following pair: $\{(e_{2|1},e_{1|2})\}$. Intuitively, the prediction mechanism learns the transition $e_{2|1}\rightarrow e_{1|2}$ thus,  when $e_{2|2}$ is available it should predict $e_{1|3}$. Note further that the sum of the subscripts addresses the same point in "simulation-prediction" time as $(k+i)$.  

The most complex structure of the prediction framework suggests a training process of $(n_x\cdot N)$ different GPs that might be challenging to train online. However, several improvements can be suggested in making it faster, as in \cite{Crocetti2023} or by adapting only the hyper-parameters of the posterior. Moreover, the training data can scale in the past up to the horizon length, which will accelerate the process further.  
\subsection{The LPVMPC-GPs framework}
The first part of the controller will stabilize the closed loop, and the MPC design will drive the system to a given reference. The whole control problem can be cast as a constrained optimization problem within a given prediction horizon of length $N$. The energy (cost) $J_k$ can be penalized with the quadratic weighted matrices\footnote{The quadratic weighted cost is defined as $\lVert x\rVert_Q^2=x^{\top}Qx$. Similarly, for the weighted norms with $R$ and $P$.}. The quadratic costs for the state and input within the prediction horizon are $Q,~R$, and the terminal cost is $P$ as the solution of the Lyapunov equation in \cite{PANDEY20178624}. At $t_k=k\cdot t_s$ and for $i=0,\ldots,N-1$, the optimization problem results in a classical quadratic program (QP) with efficient algorithms that reach real-time performance. The energy function for the regulation problem i.e., $x^{\text{ref}}=0$, to be minimized is
\begin{equation}
\footnotesize
\begin{aligned}
    J_k(u_{i|k})&:=\underset{u_{i|k}^*}{\min}\sum_{i=0}^{N-1}\left(\lVert \hat{x}_{i|k}\rVert_Q^2+\lVert u_{i|k}\rVert_R^2\right)+\lVert \hat{x}_{N|k}\rVert_P^2.
\end{aligned}
\end{equation}
In addition, within the QP problem, we can introduce linear constraints. These consist of the following sets:
\begin{equation}\label{eq:constraints}
    \begin{aligned}
        &\hat{x}_{i|k}\in\mathcal{X}_{i|k}=\{\hat{x}_k\in\IR^n|G_k^x \hat{x}_k\leq h_k^x\},\\
        &u_{i|k}\in\mathcal{U}_{i|k}=\{u_k\in\IR^{m}|G_k^u u_k\leq h_k^u\}.\\
    \end{aligned}
\end{equation}
\begin{problem}{QP function as $\texttt{QP}(\hat{p}_{i|k},{x}_k,x_{i|k}^{\text{ref}},\hat{e}_{i+1|k})$}\label{prob:Prob_LPVMPC}
\begin{subequations}\label{eq:LPV_MPC}
\footnotesize
	\begin{align} 
		\underset{u_{i|k}^*}{\text{min}} \
		& \! \sum_{i=0}^{N-1}\left(\lVert \bar{x}_{i|k} - x^{\text{ref}}_{i|k}\rVert ^2_Q + \lVert u_{i|k}\rVert^2_R\right) + \lVert \bar{x}_{N|k} - x^{\text{ref}}_{N|k}\rVert ^2_P   \\
		 \text{s.t.} \;\;
		& \hat{x}_{i+1|k}\!=\!\!A_c(\hat{p}_{i|k})\hat{x}_{i|k}\!\!+\!\!Bu_{i|k},~i\!=\!0,\!\ldots\!,\!N\!\!-\!\!1 \\ \label{eq:errorxh}
        & \bar{x}_{i+1|k}=\hat{x}_{i+1|k}+\hat{e}_{i+1|k} \\
		& \hat{x}_{0|k} = x_{0|k}=x_k, \\
		& \bar{x}_{i|k} \in \mathcal{X}_{i|k}\ominus\mathbb{E}_{i|k},~\text{with the probabilistic set}~\mathbb{E}_{i|k}\ni e_{i|k}  \label{eq:NMPC_state_cons}\\
		& u_{i|k} \in \mathcal{U}_{i|k} \label{eq:NMPC_input_cons}, \quad \forall i = 0,1,\ldots,N-1.
	\end{align}
\end{subequations}
\end{problem}

The decision variables in \cref{prob:Prob_LPVMPC} can be considered explicitly the control input and implicitly the states. The state and input constraints in~\cref{eq:NMPC_state_cons} and \cref{eq:NMPC_input_cons} are defined in~\cref{eq:constraints}. The predictor LPV model \eqref{eq:errorxh} has been updated with an error that will correct the deviation between the estimated and true response of the system. To solve the optimization problem labeled as \cref{prob:Prob_LPVMPC} optimally, the estimated scheduling signal $\hat{p}_{i|k},~i=0,\ldots,N-1$ needs to be numerically substituted and provided as a prediction, which will inevitably introduce some error. The prediction of the scheduling signal at time step $k$ is clarified. The concept is to estimate at time $k$ using the previous true (simulated) response of the actual plant (equivalent LPV that can be simulated quickly) at time $k-1$ by ensuring $\hat{p}_{i|k}=p_{i+1|k-1}=\rho(x_{i+1|k-1})$. This leverages information from the previous true response of the system and yields a good prediction due to the assumed smoothness (smooth switching) of the involved nonlinear operators. To enhance this information, we predict the forward error $e_{i+1|k}$ using GPs, which will improve the state estimation to be closer to the true one. \Cref{alg:LPVMPC} concisely summarizes the computational procedure.
\begin{algorithm}
\caption{The QP-based LPVMPC algorithm}
\textbf{Input}: Initial conditions $x_k$, the reference $(x^{\texttt{ref}},y^{\texttt{ref}})$\\
\textbf{Output}: The control input $u_k,~k=0,1,\ldots$, that drives the nonlinear system to the reference under linear constraints.
\begin{algorithmic}[1]\label{alg:LPVMPC}
\State Initialize for $k$ the scheduling parameter $\hat{p}_{i|k}$ as
$$\hat{p}_{i|k}:=\rho\left(x_k\right),~i=0,\ldots,N-1$$ 
\While{$k=0,1,\ldots$} 
\State Update the state reference $x_{i|k}^{\text{ref}}$ 
\State Learn GPs as in \eqref{eq:GPlearnPredict} and predict $\hat{e}_{i+1|k}$ 
\State Solve the QP in \eqref{eq:LPV_MPC} for $i=0,\ldots,N-1$
\begin{equation*}
\begin{aligned}    
\left[\bar{x}_{i+1|k},u_{i|k}\right]&\leftarrow\texttt{QP}^{(j)}(\hat{p}_{i|k},x_k,x_{i|k}^{\text{ref}},\hat{e}_{i+1|k}),\\
\text{Update}~\hat{p}_{i|k}&:=\rho(x_{i|k}),~i=0,\ldots,N
\end{aligned}
\end{equation*}
  \State Apply $u_k=u_{0|k}$ to the system \eqref{sys:LPVLQR}
  \State $k\leftarrow k+1$
  \State Measure $x_{k}$
  \State Update $\hat{p}_{i|k}=\hat{p}_{i+1|k-1},~i=0,\ldots,N-1$
  \EndWhile
\end{algorithmic}
\end{algorithm}
\section{Results}\label{sec:results}
\textbf{The unbalanced disk regulation problem:} We start by providing the \cref{tab:MPCPar} with all the input-output constraints and tuning parameters for solving the LPVMPC problem.
\begin{table}[h]
\footnotesize
    \centering
    \caption{MPC Parameters}
    \label{tab:MPCPar}
   \setlength{\tabcolsep}{2pt}
    \begin{tabular}{l r | l r} 
        \textbf{Parameter}   & \textbf{Value}  & \textbf{Parameter}    & \textbf{Value}\\ 
                       Lower bound on  $\theta_k$ &   $-2\pi$ [rad] & Upper bound on  $\theta_k$ &   $2\pi$ [rad]  \\
        Lower bound on  $\omega_k$ &   $-10\pi$[rad/s] & Upper bound on  $\omega_k$ &   $10\pi$[rad/s] \\
        Lower bound on  $u_k$ &   $-10$ [V] & Upper bound on  $u_k$ &  $10$ [V]  \\
           Sampling time $t_s$ & $0.01$ [s] &   Horizon length $N$ & $10$    \\
           Quad. state cost $Q$ &$\texttt{diag}(8,0.1)$ & Quad. input cost & $R=0.5$\\
    \end{tabular}
\end{table}
The terminal cost $P$ is computed from the solution of the Lyapunov matrix equation in the robust model-based LQR together with the feedback gain $K$ from \cite{PANDEY20178624}. The unbalanced disk control regulator problem is under consideration that starts from the initial conditions $x_0=\left[\begin{array}{cc}
    -2\pi & 0 \\
\end{array}\right]^\top$. We want to drive the dynamics to the state origin that is our reference signal $x^{\text{ref}}=\left[\begin{array}{cc}
   0  & 0 \\
\end{array}\right]^\top$. The dynamical system that describes the phenomenon with angular displacement-$\theta(t)$ and angular speed-$\omega(t)$ (i.e., $\omega(t)=\dot{\theta}(t)$) with state vector $x(t)=[\begin{array}{cc}
   \theta(t)  & \omega(t)   
\end{array}]^\top$ has the continuous space-state representation as in \eqref{sys:LPV} after introducing the scheduling variable $p(t):=\sin(\theta(t))/\theta(t):=\sinc(\theta(t))$. The matrices that define the continuous in-time system with $p(t)=\sinc(\theta(t))$ are
\begin{equation}
\footnotesize
    A_{cont}(p(t))=\left[\begin{array}{cc}
        0 & 1 \\
       \frac{mgl}{I_n}p(t)  & -\frac{1}{\tau}
    \end{array}\right],~B_{cont}=\left[\begin{array}{c}
         0  \\
          \frac{K_m}{\tau}
    \end{array}\right],
\end{equation}
and the chosen parameters are provided in \cref{tab:disk_parameters}.
\begin{table}
    \centering
     \caption{Parameters for the unbalanced disk example}
    \begin{tabular}{c|l}
        $I_n$ &  $2.4\cdot1e-4$ [kg$\cdot m^2$]\\
        $m$ &  $0.076~$ [Kg]\\
        $g$ & $9.81$ [m/s]\\
        $l$ & $0.041$ [m]\\
        $\tau$ & $0.4$ [1/s] \\
        $K_m$ & $11$ [rad/Vs$^2$]
    \end{tabular}
    \label{tab:disk_parameters}
\end{table}
We discretize with forward Euler\footnote{Forward Euler: $\dot{x}(t_k)\approx\frac{x(t_k+t_s)-x(t_k)}{t_s},~t_k=t_s\cdot k,~k\in\IZ_+$.}, and with $x_k=\left[\begin{matrix}
   \theta_k & \omega_k \end{matrix}\right]^\top$, the remaining discrete LPV system \eqref{sys:LPV} with scheduling parameter $p_k=\sinc(\theta_k)$ has the following matrices
\begin{equation}
\footnotesize
\begin{aligned}\label{eq:discdisc}
    A(p_k)&=\left[\begin{array}{cc}
        1 & t_s \\
       t_s\frac{mgl}{I_n}p_k  & 1-\frac{t_s}{\tau}
    \end{array}\right],~B=\left[\begin{array}{c}
         0  \\
          t_s\frac{K_m}{\tau}
    \end{array}\right].
    \end{aligned}
\end{equation}
\vspace{5mm}
\begin{figure*}[t]
  \centering
  \includegraphics[scale=0.27]{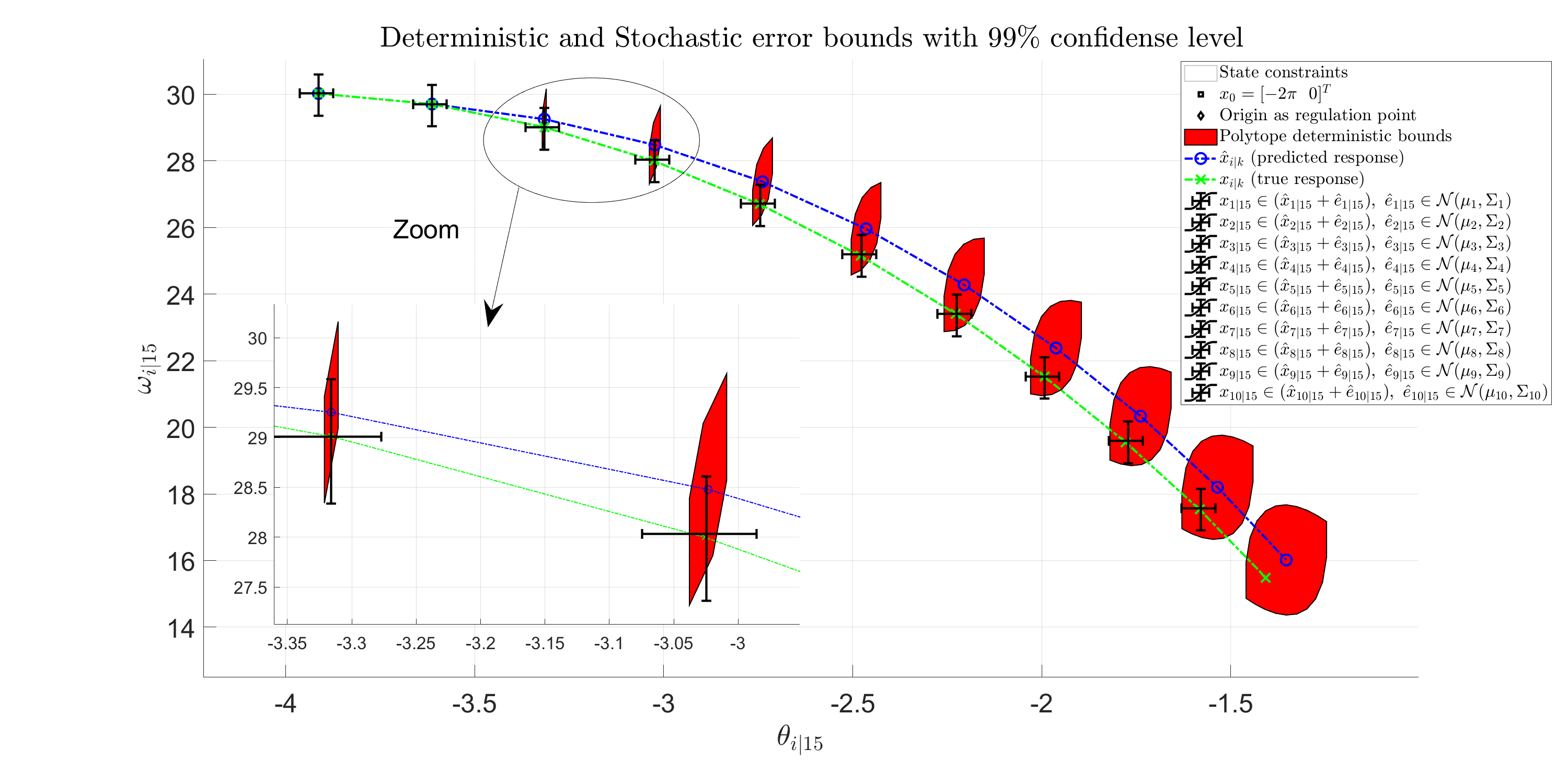}
     \caption{Phase state space evolution of the regulation problem starting from $x_0=[-2\pi,0]^\top$ with the true (green) and predicted (blue) responses of the system at time $k=15$ and for the horizon $i=0,\ldots,10$. Deterministic error polytopes (red) centered at predicted values enclosing the true response, and the error bars indicate a probabilistic region (2D Gaussian distribution) computed from the GPs uncertainty quantification framework and successfully reduces conservatism.}
    \label{fig:errorbounds}
\end{figure*}
In \cref{fig:errorbounds}, the stochastic bounds that can be computed through the GPs regression framework succeeded in reducing the conservatism and complexity (Minkowski sums) of the deterministic (red polytopes) successfully from \cite{karachalios2023error}. We know that the error will not be out of the red polytopes in \cref{fig:errorbounds}; the error bars indicate a $99\%$ confidence interval that, indeed, the error stays within. At time $k=15\cdot t_s=0.15(s)$, the error $e_{i|15}$ is quantified for the future prediction points $i=1,\ldots,N=10$, with a mean, variance and for both states $(\theta,\omega)$ that results to a sequence of joint Gaussian distributions. The estimation of the true response is quite accurate $\bar{x}\approx x$, and with high confidence, the variability is much less than in the determinist approach, capable of handling better the system's disturbances, certifying robustness.

\begin{figure}[!ht]
    \centering
    \includegraphics[scale=0.16]{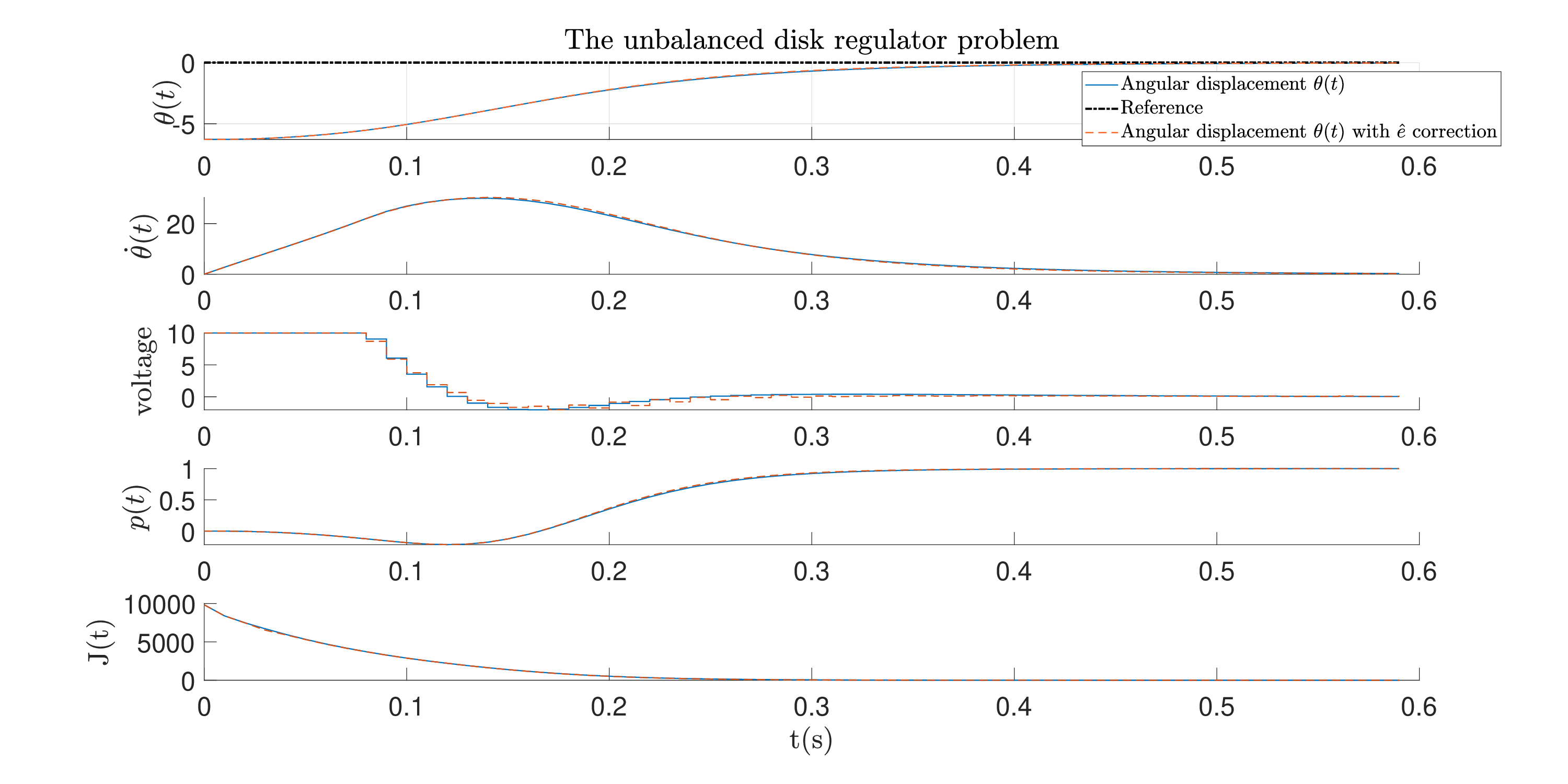}
    \caption{State, control, and scheduling trajectories of the regulation problem along with the optimization cost for the two cases with (orange) or without (cyan) error consideration in the LPVMPC solution algorithm. The regulator problem solved with a fair approximation at $k=60$ which translates to real-time $0.6$ (s).}
    \label{fig:fig3}
\end{figure}

\begin{figure}[!ht]
    \centering
    \includegraphics[scale=0.16]{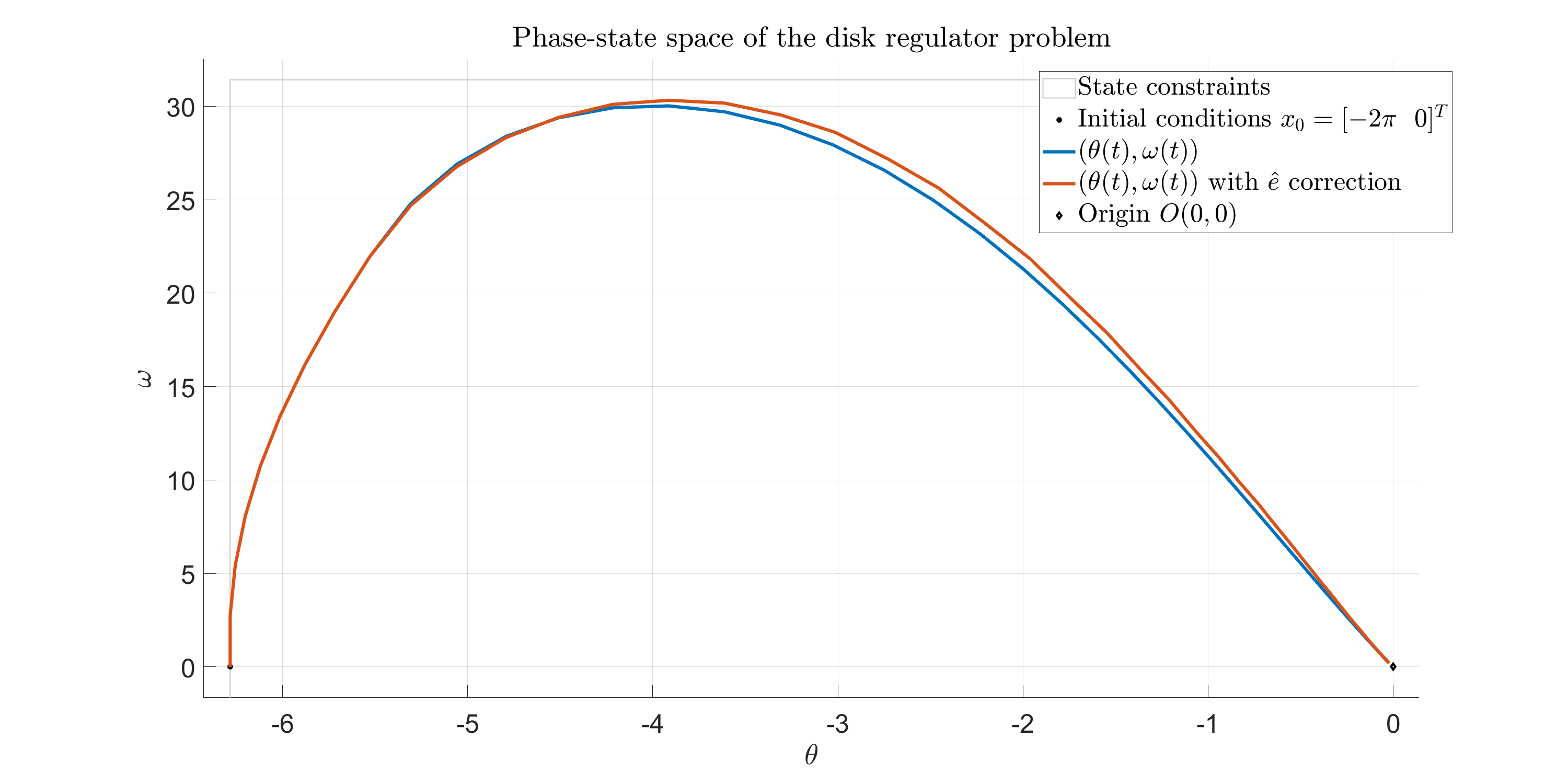}
    \caption{Phase-state trajectory evolution with a slight improvement in the performance has been achieved when considering the predicted error from the GPs-based framework within the LPVMPC.}
    \label{fig:PhaseSpace}
\end{figure}

In \cref{fig:fig3}, both solutions of the LPVMPC problem with and without error correction for the regulator problem of the unbalanced disk are illustrated and compared. When the error estimation with its bounds is considered, we didn't notice performance loss due to the robustification using the error bounds \cref{fig:PhaseSpace}. However, we can justify the error bounds and use this information for the robustification of the MPC scheme. The monotonicity of the angular displacement $\theta(t)$ that reaches the origin target without overshooting outlines the good performance seen in NMPC frameworks. In addition, the computational burden has been avoided after utilizing the QP performance, and the GP framework is promising for handling disturbance sets and faster computations compared to the deterministic bound analysis.   

\section{CONCLUSIONS}\label{sec:conclusion}
In this study, we predicted the error propagation between the true and predicted responses with Gaussian Processes (GPs) by providing stochastic bounds that successfully reduce conservatism compared to a recent deterministic approach based on polytopic tube bounds. The derived method that combines the LPVMPC and GPs frameworks solves nonlinear control problems by offering variance measures that will certify robustness against the system's disturbances and/or measurement noise. 

The study on reducing the number of trained GPs by exploiting the error mechanism that might inherit some global regularity has been left for our future research endeavors. Together with that being mentioned, we are targeting to solve the error prediction problem online and provide constraint tightening that could improve performance and reduce conservatism. Finally, developing such tools, even in the stochastic sense, will support our long-term goal of providing theoretical guarantees such as stability and recursive feasibility to assert safety in autonomous systems that can vary from mechanical to medical engineering applications.

\bibliographystyle{IEEEtran}
\bibliography{mybib}
\end{document}